\documentclass[11pt]{article}
\usepackage{amsfonts,mathrsfs,amssymb,amsthm,mathptm}
\usepackage{amsmath,amscd}
\usepackage{mathptm,pslatex}
\usepackage{color}
\usepackage[dvips]{graphicx}
\usepackage[all]{xy}
\usepackage{graphicx}
\usepackage{fancyhdr}
\usepackage{url}

\begin{document}
%\sloppy
\newtheorem{Def}{Definition}[section]
\newtheorem{Bsp}[Def]{Example}
\newtheorem{Prop}[Def]{Proposition}
\newtheorem{Theo}[Def]{Theorem}
\newtheorem{Lem}[Def]{Lemma}
\newtheorem{Koro}[Def]{Corollary}
\theoremstyle{definition}
\newtheorem{Rem}[Def]{Remark}

\newcommand{\add}{{\rm add}}
\newcommand{\con}{{\rm con}}
\newcommand{\gd}{{\rm gl.dim}}
\newcommand{\dm}{{\rm domdim}}
\newcommand{\tdim}{{\rm dim}}
\newcommand{\E}{{\rm E}}
\newcommand{\Mor}{{\rm Morph}}
\newcommand{\End}{{\rm End}}
\newcommand{\ind}{{\rm ind}}
\newcommand{\rsd}{{\rm res.dim}}
\newcommand{\rd} {{\rm rep.dim}}
\newcommand{\ol}{\overline}
\newcommand{\overpr}{$\hfill\square$}
\newcommand{\rad}{{\rm rad}}
\newcommand{\soc}{{\rm soc}}
\renewcommand{\top}{{\rm top}}
\newcommand{\stp}{{\mbox{\rm -stp}}}
\newcommand{\pd}{{\rm projdim}}
\newcommand{\id}{{\rm injdim}}
\newcommand{\fld}{{\rm flatdim}}
\newcommand{\fdd}{{\rm fdomdim}}
\newcommand{\Fac}{{\rm Fac}}
\newcommand{\Gen}{{\rm Gen}}
\newcommand{\fd} {{\rm findim}}
\newcommand{\Fd} {{\rm Findim}}
\newcommand{\Pf}[1]{{\mathscr P}^{<\infty}(#1)}
\newcommand{\DTr}{{\rm DTr}}
\newcommand{\cpx}[1]{#1^{\bullet}}
\newcommand{\D}[1]{{\mathscr D}(#1)}
\newcommand{\Dz}[1]{{\mathscr D}^+(#1)}
\newcommand{\Df}[1]{{\mathscr D}^-(#1)}
\newcommand{\Db}[1]{{\mathscr D}^b(#1)}
\newcommand{\C}[1]{{\mathscr C}(#1)}
\newcommand{\Cz}[1]{{\mathscr C}^+(#1)}
\newcommand{\Cf}[1]{{\mathscr C}^-(#1)}
\newcommand{\Cb}[1]{{\mathscr C}^b(#1)}
\newcommand{\Dc}[1]{{\mathscr D}^c(#1)}
\newcommand{\K}[1]{{\mathscr K}(#1)}
\newcommand{\Kz}[1]{{\mathscr K}^+(#1)}
\newcommand{\Kf}[1]{{\mathscr  K}^-(#1)}
\newcommand{\Kb}[1]{{\mathscr K}^b(#1)}

\newcommand{\modcat}{\ensuremath{\mbox{{\rm -mod}}}}
\newcommand{\Modcat}{\ensuremath{\mbox{{\rm -Mod}}}}
\newcommand{\stmodcat}[1]{#1\mbox{{\rm -{\underline{mod}}}}}
\newcommand{\pmodcat}[1]{#1\mbox{{\rm -proj}}}
\newcommand{\imodcat}[1]{#1\mbox{{\rm -inj}}}
\newcommand{\Pmodcat}[1]{#1\mbox{{\rm -Proj}}}
\newcommand{\Imodcat}[1]{#1\mbox{{\rm -Inj}}}
\newcommand{\opp}{^{\rm op}}
\newcommand{\otimesL}{\otimes^{\rm\mathbb L}}
\newcommand{\rHom}{{\rm\mathbb R}{\rm Hom}\,}
\newcommand{\projdim}{\pd}
\newcommand{\Hom}{{\rm Hom}}
\newcommand{\Coker}{{\rm Coker}}
\newcommand{ \Ker  }{{\rm Ker}}
\newcommand{ \Cone }{{\rm Con}}
\newcommand{ \Img  }{{\rm Im}}
\newcommand{\Ext}{{\rm Ext}}
\newcommand{\StHom}{{\rm \underline{Hom}}}
\newcommand{\gm}{{\rm _{\Gamma_M}}}
\newcommand{\gmr}{{\rm _{\Gamma_M^R}}}
\def\vez{\varepsilon}\def\bz{\bigoplus}  \def\sz {\oplus}
\def\epa{\xrightarrow} \def\inja{\hookrightarrow}

\newcommand{\lra}{\longrightarrow}
\newcommand{\llra}{\longleftarrow}
\newcommand{\lraf}[1]{\stackrel{#1}{\lra}}
\newcommand{\llaf}[1]{\stackrel{#1}{\llra}}
\newcommand{\ra}{\rightarrow}
\newcommand{\dk}{{\rm dim_{_{k}}}}
\newcommand{\colim}{{\rm colim\, }}
\newcommand{\limt}{{\rm lim\, }}
\newcommand{\Add}{{\rm Add }}
\newcommand{\Tor}{{\rm Tor}}
\newcommand{\Cogen}{{\rm Cogen}}
\newcommand{\Tria}{{\rm Tria}}
\newcommand{\tria}{{\rm tria}}

{\Large \bf
\begin{center}
Frobenius bimodules and flat-dominant dimensions
\end{center}}

\medskip
\centerline{\textbf{Changchang Xi}}

\renewcommand{\thefootnote}{\alph{footnote}}
\setcounter{footnote}{-1} \footnote{ $^*$ Corresponding author.
Email: xicc@cnu.edu.cn; Fax: 0086 10 68903637.}
\renewcommand{\thefootnote}{\alph{footnote}}
\setcounter{footnote}{-1} 
\renewcommand{\thefootnote}{\alph{footnote}}
\setcounter{footnote}{-1} 

\begin{abstract}
We establish relations between Frobenius parts and between flat-dominant dimensions of algebras linked by Frobenius bimodules. This is motivated by the Nakayama conjecture and an approach of Martinez-Villa to the Auslander-Reiten conjecture on stable equivalences. We show that the Frobenius parts of Frobenius extensions are again Frobenius extensions. Further, let $A$ and $B$ be finite-dimensional algebras over a field $k$, and let $\dm(_AX)$ stand for the dominant dimension of an $A$-module $X$.  If $_BM_A$ is a Frobenius bimodule, then $\dm(A)\le \dm(_BM)$ and $\dm(B)\le \dm(_A\Hom_B(M, B))$. In particular,
if $B\subseteq A$ is a left-split (or right-split) Frobenius extension, then $\dm(A)=\dm(B)$. These results are applied to
calculate flat-dominant dimensions of a number of algebras: shew group algebras,
stably equivalent algebras, trivial extensions and Markov extensions. Finally, we prove that the universal (quantised) enveloping algebras of semisimple Lie algebras are $QF$-$3$ rings in the sense of Morita.
\end{abstract}

\section{Introduction}

Let $R$ and $S$ be rings with identity. Recall that a $S$-$R$-bimodule $M$ is called a \emph{Frobenius bimodule} if it is finitely generated
and projective as a left $S$- and right $R$-module, and there is an $R$-$S$-bimodule
isomorphism $$^*M := {}_R\Hom_S(M,S)_S\simeq {}_R\Hom_{R^{\opp}}(M_R, R_R)_S =: M^*.$$ A prominent example of Frobenius bimodules is provided by Frobenius extensions. Recall that an extension $S\subseteq R$ of rings (or algebras) is called a \emph{Frobenius extension (or a free Frobenius extension)} if $_SR$ is a finitely generated projective $S$-module (or a free $S$-module of finite rank) and $_RR_S\simeq \Hom_S(_SR,{}_SS)$ as $R$-$S$-bimodules (see \cite{Kasch1} and \cite{Nakayama2}).
This is equivalent to saying that $R_S$ is a finitely generated projective right $S$-module and $_SR_R \simeq \Hom_{S^{\opp}}(R_S,S_S) $ as $S$-$R$-bimodules.
In this case, both $_SR_R$ and $_RR_S$ are Frobenius bimodules. If $S$ is a field and $R$ is a finite-dimensional $S$-algebra, then the notion of Frobenius extensions is just the definition of Frobenius algebras. Frobenius algebras and extensions are of broad interest in many different areas (see, for instance, \cite{Kadison}), they not only give rise to Frobenius functors for constructing interesting topological
quantum field theories in dimension 2 and even 3 (for example, see \cite{37}) and provide connections
between representation theory and knot theory (for example in the spirit of \cite{Kadison1996}), but also are used in the study of Calabi-Yau properties of Cherednik algebras and quantum algebras (see \cite{BrownGordonStroppel}).

Dominant dimensions of finite-dimensional algebras were initiated by Nakayama in \cite{Nakayama}, where he suggested to classify finite-dimensional algebras according to the length of their minimal injective resolutions with the first $n$ terms being projective. They were further investigated substantially by Tachikawa, M\"uller and others (see \cite{Tachikawa, Muller}). Behind the notion of dominant dimensions is the famous Nakayama conjecture: A finite-dimensional algebra over a field with infinity dominant dimension is self-injective. This is a long-standing but still not yet solved open question (see \cite{Nakayama}, \cite[Conjecture (8), p. 410]{Auslander1995}). Recently, dominant dimensions are tied up with tilting modules (see \cite{CX6}).

In this note, we shall establish relations between (flat-)dominant dimensions and between Frobenius parts of algebras (or rings) linked by Frobenius bimodules, including Frobenius extensions of algebras.

To facilitate our consideration, let us recall some basic terminology.

By $R\Modcat$ (respectively, $R\modcat$) we denote the category of all left (respectively, finitely generated) $R$-modules. By $A$ and $B$ we denote Artin algebras over a fixed commutative Artin ring $k$. In this case we only consider finitely generated left $A$-modules for an Artin algebra $A$.

For an Artin algebra $A$ and an $A$-module $_AM$, the dominant dimension of $_AM$, denoted by $\dm(_AM)$, is by definition the minimal number $n$ in a minimal injective resolution of $_AM$: $$ 0\lra {}_AM\lra I^0\lra \cdots\lra I^n\lra \cdots $$such that $I^n$ is not projective. If such a number does not exist, we write $\dm(_AM)=\infty$. The dominant dimension of $A$ is defined to be $\dm(_AA)$, denoted $\dm(A)$. For two $A$-modules $X$ and $Y$, $\dm(X\oplus Y)=min\{\dm(X),\dm(Y)\}$.
If $A$ is a self-injective algebra, that is, the regular $A$-module $_AA$ is injective, then $\dm(A)=\infty$. The converse of this statement is precisely the foregoing mentioned Nakayama conjecture.

Given an Artin algebra $A$, let $\nu$ be the Nakayama functor $D\Hom_A(-, {}_AA)$, where $D$ is the usual duality of an Artin algebra. Recall that a projective $A$-module $P$ is said to be $\nu$-\emph{stably projective} if $\nu^iP$ is projective for all $i> 0$, where $\nu$ is the Nakayama functor $D\Hom_A(-,{}_AA)$. We denote by $A$-stp the full subcategory of $A$-mod consisting of all $\nu$-stably projective $A$-modules. Following \cite{FHK}, we denote by $\nu$-$\dm(M)$ the minimal number $n$ in the minimal injective resolution of an $A$-module $_AM$ such that
$I^n$ is not $\nu$-stably projective. Clearly, $\nu$-$\dm(M)\le \dm(M)$. As in \cite{HuXi2015}, the \emph{Frobenius part} of $A$, which is unique up to Morita equivalence, is defined as the endomorphism algebra of an $A$-module $X$ with $\add(X)= A$-stp, where $\add(X)$ stands for the additive subcategory of $A$-mod generated by $X$.  Frobenius parts of algebras play a significant role in the study of the Auslander-Reiten conjecture which says that the numbers of non-projective simple modules of stably equivalent algebras are equal (see \cite{Martinez-Villa1990a} for Martinez-Villa's approach to this conjecture), dominant dimensions (see \cite{CX6}) and lifting stable equivalences to derived equivalences (see \cite{HuXi2015}).

Generally, for an arbitrary ring $R$ and an $R$-module $_RM$, we can similarly define the dominant dimension of $M$. But we need to distinguish left and right dominant dimensions. Following \cite{Hoshino}, the \emph{left flat-dominant dimension} of $M$ is the maximal nummber $n$ such that all $I^0, \cdots, I^{n-1}$ are flat $R$-modules in a minimal injective resolution of $_RM$:
$$0\lra {}_RM\lra I^0\lra \cdots\lra I^n\lra \cdots $$
We denote by $\fdd(M)$ the left flat-dominant dimension of $M$. The \emph{left flat-dominant dimension} of $R$ is defined to be $\fdd(_RR)$ and denoted by $\fdd(R)$. If $R$ is an Artin algebra, then $\fdd(R)=\dm(R)$. If $R$ is a noetherian ring (that is, a left and right noetherian ring), then the left and right flat-dominant dimensions of $R$ are equal (see \cite{Hoshino}). For the ring $\mathbb{Z}$ of integers, $\fdd(\mathbb{Z})=1$ and $\dm(\mathbb{Z})=0$. This follows from the injective resolution $0\ra \mathbb{Z}\ra \mathbb{Q}\ra \mathbb{Q}/\mathbb{Z}\ra 0$ and the fact that a $\mathbb{Z}$-module is flat if and only if it is torsion-free.

Flat-dominant dimensions were used in \cite{T} to study almost coherence and left FTF property (see \cite[Section 1]{T} for definition) for (quasi-) Frobenius extensions $S\subseteq R$ of rings. For instance, it was proved in \cite{T} that the ring $S$ is left FTF if and only if $R$ is a left FTF.

By a \emph{flat-injective} module we mean a module which is both flat and injective. An extension $S\subseteq R$ of rings is said to be \emph{left-split} (or \emph{right-split}) if the inclusion map is a split monomorphism of left (or right) $S$-modules. For example, if $S\subseteq R$ is a Frobenius extension of rings with $S$ a commutative ring, then the extension is both left- and right-split (see \cite[III.4.8, Lemma 2]{browngoodearl}).

We first state a result on Frobenius parts of algebras linked by a Frobenius bimodule.

\begin{Theo}\label{frobpart}
Let $_BM_A$ be a Frobenius bimodule of Artin algebras $A$ and $B$ such that both $_BM$ and $M_A$ are faithful, and let $X\in A\modcat$ and $Y\in B\modcat$ such that $\add(X)=A\emph{-stp}$ and $\add(Y)=B\emph{-stp}$. Then there exist canonical injective homomorphisms $\End_A(X)\ra \End_B(M\otimes_AX)£©$ and $\End_B(Y)\ra \End_A(\Hom_B(_BM_A,{}_BY))$ of endomorphism algebras such that they are both left- and right-split, Frobenius extensions.
\end{Theo}

Next, we compare the flat-dominant dimensions of algebras linked by Frobenius bimodules, including Frobenius extensions.

\begin{Theo}
$(1)$  Let $R$ and $S$ be rings and $_SM_R$ be a Frobenius bimodule. Then $\fdd(R)\le \fdd(_SM)$ and $\fdd(S)\le \fdd(_R{}^*M)$. In particular, if $S\subseteq R$ is a free Frobenius extension or a Frobenius extension with $S$ a commutative ring, then $\fdd(R)=\fdd(S)$ and $\dm(R)=\dm(S)$.

$(2)$ Let $A$ and $B$ be Artin algebras. If $_BM_A$ is a Frobenius bimodule, then $\dm(A)\le \dm(_BM)$ and $\dm(B)\le \dm(_A{}^*M)$. Further, $\nu$-$\dm(A)\le \nu$-$\dm(_BM)$ and $\nu$-$\dm(B)\le \nu$-$\dm(_A{}^*M)$.
\label{thmdm1}
\end{Theo}

By Theorem \ref{thmdm1}, if $B\subset A$ is a Frobenius extension of Artin algebras, then we consider the Frobenius bimodule $_BA_A$ and get $\dm(B)\le \dm(A)$ and $\nu$-$\dm(B)\le \nu$-$\dm(A)$. In particular, if the Frobenius extension is left-split, then $\dm(B)=\dm(A)$. Thus, for a left-split Frobenius extension $B\subseteq A$, if the Nakayama conjecture holds true for $B$, then so does it for $A$. This is due to the fact that, for a Frobenius extension $B\subseteq A$ of Artin algebras, if $B$ is a self-injective algebra, then $A$ is also a self-injective algebra.

We will apply Theorem \ref{thmdm1} to to skew group algebras, trivial extensions, Markov extensions and stably equivalent algebras as well as to the  universal enveloping algebras of semisimple Lie algebras.

The proofs of the above theorems will be carried out in the next sections.

\section{Proofs of the results}

Throughout this paper, a ring (or an algebra) means an associative ring (or algebra) with identity. We denote by $R^{\opp}$ the opposite ring of a ring $R$.
For two homomorphisms $f:X\ra Y$, $g: Y\ra Z$ of $R$-modules, the composition of $f$ with $g$ is denoted by $fg$ which is a homomorphism from $X$ to $Z$. So the image of $x\in X$ under $f$ is denoted by $(x)f$ instead of $f(x)$.
In this way, $\Hom_R(X,Y)$ admits automatically a left $\End_R(X)$- and right $\End_R(Y)$-bimodule.

By an Artin algebra we mean an algebra over a commutative Artin ring $k$ such that it is a finitely generated $k$-module.

\subsection{Proof of Theorem \ref{frobpart}}

We first display some properties of Frobenius bimodules. Suppoes that $R$ and $S$ are rings.
Let $_SM_R$ be a Frobenius bimodule, and we define $N:= {}^*M$. Then

(1) $M\otimes_R- \simeq \Hom_R(_RN,-): R\modcat\ra S\modcat$; $N\otimes_S-\simeq \Hom_S(M,-): S\modcat\ra R\modcat.$

(2) If an $R$-module $X$ is projective or injective, then so is the $S$-module $M\otimes_RX$.

(3) $\Ext_S^n(X,M\otimes_RY)\simeq \Ext^n_R(\Hom_R(M,X),Y)$ and $\Ext_R^n(Y,N\otimes_SX)\simeq \Ext^n_S(\Hom_R(N,Y),X)$ for all $n\ge 0$, $S$-modules $X$ and $R$-modules $Y.$ In particular, if $\Gamma\subseteq \Lambda$ is a Frobenius extension of rings, then $_{\Lambda}\Lambda_{\Gamma}$ is a Frobenius bimodule and there are isomorphisms of abelian groups:
$ \Ext_{\Lambda}^n(U, \Lambda\otimes_{\Gamma}V)\simeq \Ext^n_{\Gamma}(U,V)$ and $\Ext^n_{\Gamma}(V, U)\simeq \Ext^n_{\Lambda} (\Hom_{\Gamma}(\Lambda, V), U)$
for all $n\ge 0$, $\Lambda$-modules $U$ and $\Gamma$-modules $V$; and it follows that $\pd_{\Lambda}(U)=\pd_{\Gamma}(U)$ for all $\Lambda$-modules $U$ (see also \cite[(IX)]{Kasch1}).

($3'$) $\Ext^i_S(M\otimes_RX,Y)\simeq \Ext^i_R(X, N\otimes_BY)$ and $\Ext^i_R(N\otimes_SY,{}_RX)\simeq \Ext^i_S(Y, M\otimes_RX)$ for all $i\ge 0$, $_RX$ and $_SY$.

Now, we assume that $A$ and $B$ are Artin algebras and that $_BM_A$ is a Frobenius bimodule with $N:= {}^*M$.

(4) $\nu_B(M\otimes_AX)\simeq M\otimes_A\nu_AX$ for all $X\in A\modcat$, and $\nu_A\Hom_B(M,{}_BY)\simeq \Hom_B(M,\nu_BY)$ for all $Y\in B\modcat$. Thus $\tau_B(M\otimes_AX)\simeq M\otimes_A\tau_AX$ for all $X\in \modcat{A}$ and $\tau_A\Hom_B(M,Y)\simeq \Hom_B(M,\tau_BY)$, where $\tau$ is the Auslander-Reiten translation.

In fact, for any $A$-module $X$, we have $$\begin{array}{rl} \nu_B(M\otimes_AX)& = D\Hom_B(M\otimes_AX, {}_BB)\\ & \simeq D\Hom_A(X, \Hom_B(M,B))\\ & = D\Hom_A(X, N) \\ & \simeq D(\Hom_A(X, A)\otimes_AN) \quad (_AN: \mbox{ finitely generated and projective }) \\ & \simeq \Hom_A(N,\nu_AX) \\ & \simeq M\otimes_A\nu_AX.\end{array}$$ The second statement of (4) can be proved similarly. The statement for $\tau$ is a consequence of the isomorphisms about $\nu$.

(5) If $M_A$ is faithful, then $A$-stp = $\add(N\otimes_BY)$ for any $_BY$ with $\add(Y) = B$-stp. Dually, if $_BM$ is faithful, then $B$-stp = $\add(M\otimes_AX)$ for any $_AX$ with $\add(X) = A$-stp.

In fact, we have $\add(N\otimes_BY)\subseteq A$-stp by (1) and (4). Now, suppose $e^2=e\in A$ such that $\add(Ae)= A$-stp. Then $M\otimes_AAe\in B$-stp by (4), that is, $Me\in \add(Y)$. Hence $\Hom_B(M, Me)\in \add(\Hom_B(M,Y))=\add(N\otimes_BY)$. Note that we have the adjunction homomorphism
$$ \mu: Ae\lra \Hom_B(M, Me), x\mapsto (\varphi_x: M\ra Me, m\mapsto mx, m\in M), x\in Ae.$$
Since $M_A$ is faithful, the above adjunction is injective. By definition, $Ae$ is also an injective $A$-module. This implies that $\mu$ is a split-monomorphism and $Ae$ is a direct summand of $N\otimes_BY$, that is, $A$-stp = $\add(Ae)\subseteq \add(N\otimes_AY)$. Thus $A$-stp = $\add(N\otimes_BY)$.

Applying (5) to Frobenius extensions, we have the following

(6) If $B\subseteq A$ is a Frobenius extension of Artin algebras, then the functor $A\otimes_B-$ induces  a functor from $B$-stp to $A$-stp, and the restriction functor induces a functor from $A$-stp to $B$-stp. Moreover, if $_AX$ is an $A$-module such that $\add(X) = A$-stp, then $B$-stp = $\add(_BX)$. Similarly, if $Y$ is a $B$-module such that $\add(Y)$ = $B$-stp, then $A$-stp = $\add(A\otimes_BY)$.

\medskip
{\bf Proof of Theorem \ref{frobpart}}:
Let $_BM_ A$ be a Frobenius bimodule such that both $_BM$ and $M_A$ are faithful. Then the Frobenius bimodule $_AN_B$ is also faithful as a left or right module. Suppose $_AX\in A\modcat$ and $_BY\in B\modcat$ with $A$-stp = $\add(X)$ and $B$-stp = $\add(Y)$. Then there are two canonical ring homomorphisms: $\varphi: \End_A(X)\ra \End_B(M\otimes_AX)$  and $\psi: \End_B(Y)\ra \End_A(N\otimes_BY)$. We show that they are injective.

To see that $\varphi$ is injective, we just note that $\varphi$ is a composition of the maps
$$ \End_A(X)\lraf{\Hom_A(X,\mu)} \Hom_A(X, \Hom_B(M,M\otimes_AX))\simeq \Hom_B(M\otimes_AX,M\otimes_AX)$$
where $\mu: X\ra \Hom_B(M,M\otimes_AX)$ is given by $x\mapsto (\alpha_x: M\ra M\otimes_AX, m\mapsto m\otimes x)$ is injective, since $X$ is of the form $Ae$ with $e^2=e$ and $M_A$ is faithful.
Similarly, we can prove that $\psi$ is injective.

(a) Define $\Gamma£º= \End_A(X)$ and $\Lambda:=\End_B(M\otimes_AX)$. Then $\Lambda$ is a $\Gamma$-bimodule whose $\Gamma$-module structures are induced by the right $\Gamma$-structure of $X_{\Gamma}$. Moreover, it follows from $\add(_BM\otimes_AX)=\add(Y)$ and $\add(N\otimes_BY)=\add(X)$ (see (5)) that
$$\begin{array}{rl}{}_{\Gamma}\Lambda & =\Hom_B(M\otimes_AX_{\Gamma},M\otimes_AX)\\
                                    & \simeq \Hom_A(_AX_{\Gamma}, \Hom_B(M,M\otimes_AX)) \\
                                    & \simeq \Hom_A(_AX_{\Gamma}, N\otimes_BM\otimes_AX)) \\
                                    \end{array}$$
is finitely generated and projective. Note that, for ${}_{\Gamma}\Lambda$ to be finitely generated and projective, one only needs $N\otimes_BM\otimes_AX\in \add(X).$ Similarly,
$$\begin{array}{rl}\Lambda_{\Gamma} & =\Hom_B(M\otimes_AX,M\otimes_AX_{\Gamma})\\
                                    & \simeq \Hom_B(M\otimes_AX, \Hom_A(_AN,X_{\Gamma})) \quad (\mbox{ by (1) })\\
                                    & \simeq \Hom_A(N\otimes_BM\otimes_AX,{}_AX_{\Gamma})\end{array}$$
is finitely generated and projective. Moreover, there are the following isomorphisms of $\Gamma$-$\Lambda$-bimodules:
$$  \begin{array}{rl}\Hom_{\Gamma^{\opp}}(_{\Lambda}\Lambda_{\Gamma},{}_{\Gamma}\Gamma_{\Gamma}) & =\Hom_{\Gamma^{\opp}}\big(\Hom_B(_BM\otimes_AX,M\otimes_AX),\Hom_A(X,X)\big)\\
                                    & \simeq \Hom_{\Gamma^{\opp}}\big(\Hom_B(_BM\otimes_AX,\Hom_A(N,X)),\Hom_A(X,X)\big)\quad ( \mbox{ by (1) })\\
                                    & \simeq \Hom_{\Gamma^{\opp}}\big(\Hom_A(N\otimes_BM\otimes_AX, {}_AX),\Hom_A(X,X)\big)\\
                                    &\simeq \Hom_A(_AX,N\otimes_BM\otimes_AX) \quad ( \mbox{ by } \Hom_A(-,X))\\
                                    & \simeq {}_{\Gamma}\Lambda_{\Lambda}.\end{array} $$

(b) Define $D:=\End_B(Y)$ and $C:=\End_A(N\otimes_BY)$. Then $C$ is a $D$-bimodule. Similarly, we can show that ${}_DC$ and $C_D$ are finitely generated and projective. Moreover, $\Hom_D(_DC,{}_DD_D)\simeq {}_CC_D.$ This is due to the isomorphisms of $C$-$D$-modules:
$$\begin{array}{rl}\Hom_D(_DC_C,{}_DD_D) & =\Hom_{D}\big(\Hom_A(_AN\otimes_BY_D, {}_A(N\otimes_BY)_C),\Hom_B(_BY_D, {}_BY_D)\big)\\
                                        & \simeq \Hom_{D}\big(\Hom_B(_BY_D,\Hom_A(_AN_B,{}_A(N\otimes_BY)_C)),\Hom_B(_BY_D, {}_BY_D)\big)\\
                                    & \simeq \Hom_{D}\big(\Hom_B(_BY_D,{}_BM\otimes_A(N\otimes_BY)_C),\Hom_B(_BY_D, {}_BY_D)\big)\quad ( \mbox{ by (1) })\\
                                    & \simeq \Hom_{B}\big(_BM\otimes_A(N\otimes_BY)_C, {}_BY_D\big) \quad ( \mbox{ by } \Hom_B(Y,-))\\
                                    &\simeq \Hom_A\big(_A(N\otimes_BY)_C, {}_AN\otimes_BY_D\big)\quad ( \mbox{ by (1) })\\
                                    & \simeq {}_{C}C_{D}.\end{array} $$

Thus $\varphi$ and $\psi$ are Frobenius extensions. Since $\End_A(X)$ and $\End_B(Y)$ both are self-injective algebras by a result of Martinez-Villa (for example, see \cite[Lemma 2.7(3)]{HuXi2015}), the two extensions $\varphi$ and $\psi$ are both left- and right-split.
$\square$

The following corollary shows that Frobenius parts of Frobenius extensions are again Frobenius extensions.
\begin{Koro}\label{frobext}
Let $B\subseteq A$ be a Frobenius extension of Artin algebras, and let $X\in A\modcat$ and $_B Y\in B\modcat$ such that $\add(X)=A\emph{-stp}$ and $\add(Y)=B\emph{-stp}$. Then the canonical injective homomorphisms $\End_A(X)\ra \End_B(X)$ and $\End_B(Y)\ra \End_A(A\otimes_BY)$ of endomorphism algebras are Frobenius extensions, respectively.
\end{Koro}

{\it Proof.} For any extension $B\subseteq A$ of algebras, the $B$-modules $_BA$ and $A_B$ are always faithful. So Corollary \ref{frobext} follows from Theorem \ref{frobpart} immediately. $\square$

\smallskip
Note that if $A$ and $B$ are $k$-algebras over a field $k$ with $B$ a Frobenius algebra (that is, ${}_BB\simeq D(B_B)$) then $A\ra A\otimes_kB$, $a\mapsto a\otimes 1$, is a Frobenius extension, where $A\otimes_kB$ is the tensor product of $A$ and $B$ over $k$. Motivated by this fact, we now calculate the Frobenius part of the tensor product of two finite-dimensional algebras over a field.

Let $A$ and $B$ be finite-dimensional $k$-algebras over a field $k$ and $\Lambda:=A\otimes_kB$ be the tensor product of $A$ and $B$ over $k$. If $_AX$ and $_BY$ are finite-dimensional modules over $A$ and $B$, respectively, then $X\otimes_kY$ is clearly an $\Lambda$-module and $D(X\otimes_kY)\simeq DX\otimes_kDY$ as right $\Lambda$-modules. This can be seen from the isomorphisms of $\Lambda$-modules:
$$\Hom_k(X\otimes_kY,k)\simeq \Hom_k(X,\Hom_k(Y,k)\simeq \Hom_k(X,k\otimes_k\Hom_k(Y,k)\simeq \Hom_k(X,k)\otimes_k\Hom_k(Y,k).$$%=DX\otimes_kDY
Thus $\nu_{\Lambda}(X\otimes_kY)=D\Hom_{\Lambda}(X\otimes_kY,A\otimes_kB)\simeq D\big(\Hom_A(X,A)\otimes_k\Hom_B(Y,B)\big)\simeq \nu_AX\otimes_k\nu_BY$, where the first isomorphism is due to \cite[XI, Theorem 3.1, pp. 209-210]{CE}. This implies that ($A$-stp)$\otimes_k(B$-stp$)\subseteq \Lambda$-stp.

Suppose that $P\in \Lambda$-stp is indecomposable. Then $P$ is of the form $X\otimes_kY$ with $X$ an indecomposable projective $A$-module and $Y$ an indecomposable projective $B$-module. Since $\nu^i_{\Lambda}P\simeq \nu_A^iX\otimes_k\nu_B^iY$, both $\nu_A^iX$ and $\nu_B^iY$ are projective for all $i\ge 0$. Thus $\Lambda\stp =(A\stp)\otimes_k(B\stp)$.

\begin{Prop} Let $A$ and $B$ be finite-dimensional $k$-algebras over a field $k$, and let $\Delta(A)$ denote the Frobenius part of $A$. Then $\Delta(A\otimes_kB)$ and $\Delta(A)\otimes_k\Delta(B)$ are Morita equivalent. \label{stp-tensorprod}
\end{Prop}

{\it Proof.} Let $_AX$ and $_BY$ be modules such that $A\stp=\add(X)$ and $B\stp=\add(Y)$. Then $\Delta(A\otimes_kB)$ is Morita equivalent to $\End_{A\otimes_kB}(X\otimes_kY)$ by the fact $\Lambda\stp =(A\stp)\otimes_k(B\stp)$. It follows from \cite[XI, Theorem 3.1, pp. 209-210]{CE} that
$\End_{A\otimes_kB}(X\otimes_kY)\simeq \End_A(X)\otimes_k\End_B(Y).$ $\square$

\smallskip
Recall that $A$ is \emph{Frobenius-free} if the Frobenius part of $A$ is zero. In many cases, Frobenius parts of algebras are non-trivial. Following \cite{HuXi2015}, we define the \emph{Frobenius type} of an algebra to be the representation type of its Frobenius part. From Proposition \ref{stp-tensorprod} we can give another proof of the fact that Frobenius types of algebras may change under derived equivalences. Let us explain this precisely. There exist finite-dimensional algebras $A$ and $B$ (see \cite[Example 5.6]{HuXi2017}), satisfying the properties:

(a) $A$ and $B$ are derived equivalent, and

(b) $A$ is Frobenius-free and $B$ is Frobenius-finite (the Frobenius part of $B$ is $C:=k[X]/(X^2)$).

By a result of Rickard, which says that derived equivalences preserve tensor products of $k$-algebras (see \cite{Rickard}), we know that $A\otimes_kC$ and $B\otimes_kC$ are derived equivalent. By Proposition \ref{stp-tensorprod}, $A\otimes_kC$ is Frobenius-free and $B\otimes_kC$  has the Frobenius part $C\otimes_kC$. Thus $B\otimes_kC$ is not Frobenius-finite, because from (b) we can say that $B\otimes_kC$ is Frobenius-tame.
So the Frobenius type of the two algebras is changed, though they are derived equivalent.
%$e\in B$, $eBe\subseteq eAe$ Frob.ext. ?

\subsection{Proof of Theorem \ref{thmdm1}}

First, we show a result on dominant dimensions for left-split or right-split extensions of rings.

\begin{Lem} Let $\Gamma\subseteq \Lambda$ be a left-split extension of rings, or a right-split extensions between noetherian rings.

 $(1)$ If both $_{\Gamma}\Lambda$ and $\Lambda_{\Gamma}$ are projective, then $\dm(\Gamma)\ge \dm(\Lambda)$.

 $(2)$ If both $_{\Gamma}\Lambda$ and $\Lambda_{\Gamma}$ are flat, then $\fdd(\Gamma)\ge \fdd(\Lambda)$.
 \label{splitext}
\end{Lem}

{\it Proof.} Note that the following facts are true for an arbitrary extension $S\subseteq R$ of rings.

(a) If $_SR$ is projective, then every projective $R$-module is also projective as $S$-module since projective $R$-module are direct summand of copies of $_RR$.

(b) If $R_S$ is projective (or flat), then every injective $R$-module is also injective as $S$-module. Indeed, for an injective $R$-module $I$, the functor $\Hom_S(-, {}_SI)=\Hom_S(-, \Hom_R(_RR_S, {}_RI))\simeq \Hom_R(R\otimes_S-, {}_RI)= \Hom_R(-, {}_RI)\circ (R\otimes_S-)$ is exact. This implies that $_SI$ is injective.

Now, assume that the given extension is left-split. For a minimal injective resolution of  a $\Lambda$-module $M$,
$$ 0\lra {}_{\Lambda}M\lra E^0\lra \cdots \lra E^n\lra \cdots, $$with $\dm(_{\Lambda}M)=n$ (or $\fdd(_{\Lambda}M=n$), we consider the restriction of this resolution to $\Gamma$-modules. Since both $_{\Gamma}\Lambda$ and ${\Lambda}_{\Gamma}$ are projective (or flat), it follows from the foregoing facts that $\dm(_{\Gamma}M)\ge n = \dm(_{\Lambda}M)$ (or $\fdd(_{\Gamma}M)\ge n = \fdd(_{\Lambda}M)$). In particular, $\dm(_{\Gamma}\Lambda)\ge \dm(\Lambda)$ (or $\fdd(_{\Gamma}\Lambda)\ge \fdd(\Lambda)$). Since the extension is left-split, $_{\Gamma}\Gamma$ is a direct summand of $_{\Gamma}\Lambda$. As a result, $\dm(\Gamma)\ge \dm(\Lambda)$ (or $\fdd(\Gamma)\ge \fdd(\Lambda)$).

The above proof works for right-split extensions and dominant dimensions defined by right modules. But since $\dm(\Lambda)=\dm(\Lambda^{\opp})$ and $\fdd(\Lambda)=\fdd(\Lambda^{\opp})$ for noetherian rings (see \cite{Hoshino}, \cite{Muller}), (1) and (2) follows. $\square$

\medskip
{\bf Proof of Theorem \ref{thmdm1}:}

Note that if $_SM_R$ is a Frobenius bimodule, then so is the $R$-$S$-bimodule $^*M$. Let $N:={}^*M$.

(1) By the definition of Frobenius bimodules, it follows that $(M\otimes_R-, N\otimes_S-)$ and $(N\otimes_S-, M\otimes_R-)$ are adjoint pairs of functors between module categories. Moreover, if $I$ is an injective $R$-module, then $_SM\otimes_RI$ is an injective $S$-module because $\Hom_S(-, {}_SM\otimes_RI)\simeq \Hom_R(N\otimes_S-, I)= \Hom_R(?,I)\circ (N\otimes_S-)$ is an exact functor. Consequently, if $I$ is a flat-injective $R$-module, then $_SM\otimes_RI$ is a flat-injective $S$-module since we have the general fact that $_SX_R\otimes_RY$ is a flat $S$-module whenever both $_SX$ and $_RY$ are flat modules.

Now, given a minimal injective resolution of $_RR$ with $\fdd(R)=n$:
$$ 0\lra {}_RR\lra I^0\lra I^1\lra I^2\lra\cdots$$
we have an exact sequence
$$ (**) \quad 0\lra {}_SM\lra M\otimes_RI^0\lra M\otimes_RI^1\lra M\otimes_RI^2\lra\cdots$$
with $_SM\otimes_RI^j$ flat-injective for all $0\le j < n$. Thus $\fdd(_SM)\ge n=\fdd(R)$. The last inequality in (1) follows from the Frobenius bimodule $_RN_S$.

Suppose that $S\subseteq R$ is a free Frobenius extension, or a left-split Frobenius extension. Then $_SS$ is a direct summand of $_SR$, %that is, $_SR$ is a free $S$-module of finite rank
and $_SR_R$ is a Frobenius $S$-$R$-bimodule. Thus $\fdd(R)\le \fdd(_SR)\le \fdd(S)$ by the first part of Theorem \ref{thmdm1}(1). Dually, $\fdd(S)\le \fdd(_R{}(^*R))=\fdd(_RR)$. Hence $\fdd(R)$ = $\fdd(S)$. Similarly, for dominant dimensions, we have $\dm(R)$ = $\dm(S)$.

In case $S$ is a commutative ring, then the Frobenius extension is left-split (see \cite[III.4.8, Lemma 2]{browngoodearl}). Thus $\fdd(R)=\fdd(S)$ and $\dm(R)$ = $\dm(S)$.

(2) This is a consequence of (1) since for Artin algebras $A$ and $M\in A\modcat$ we always have $\fdd(_AM)$ = $\dm(_AM)$. Similarly, the inequalities hold true for $\nu$-dominant dimensions. $\square$

\medskip
In the rest of this section, we consider only finite-dimensional algebras over a fixed field $k$. So, by algebras we always mean finite-dimensional algebras over $k$, and by modules always mean finite-dimensional left modules.

Recall that an algebra $A$ is called a \emph{Morita algebra} in the sense of Kerner-Yamagata (see \cite{KY}) if it is of the form $\End_{\Lambda}(X)$ where $\Lambda$ is a self-injective algebra and $_{\Lambda}X$ is a generator for $\Lambda\modcat$, that is, ${}_{\Lambda}\Lambda\in \add(X)$.

The following corollary shows that the inequality in Theorem \ref{thmdm1}(2) is optimal.

\begin{Koro} $(1)$ Suppose that $B\subseteq A$ is a Frobenius extension of $k$-algebras. If $B$ is a Morita algebra, then so is $A$.

$(2)$ If $B$ is a self-injective $k$-algebra, then $\dm(A\otimes_kB)=\dm(A)$ for any $k$-algebra $A$.\label{dmfortensor}
\end{Koro}

{\it Proof.} (1) An algebra $A$ is Morita algebra if and only if $\nu$-$\dm(A)\ge 2$ (see \cite[Proposition 2.9]{FHK}). Thus (1) follows from Theorem \ref{thmdm1}(2) since $_BA_A$ is a Frobenius bimodule for each Frobenius extension $B\subseteq A$.

(2) For any self-injective algebra $B$, we can find a Frobenius algebra $C$, that is $_CD(C)\simeq {}_CC$ as $C$-modules, such that $B$ and $C$ are Morita equivalent. Actually, $C$ is the so-called basic algebra of $B$. In this case, $A\otimes_k B$ and $A\otimes_kC$ are Morita equivalent. So we may assume that $B$ itself is a Frobenius algebra. Then $A\otimes_kB$ is a Frobenius extension of $A$. Generally, the injective map $A\ra A\otimes_kB$, $a\mapsto a\otimes 1$ for $a\in A$, is a Frobenius (symmetric, self-injective) extension  for any Frobenius (symmetric, self-injective) $k$-algebra $B$ by \cite[Lemma 4]{Sugano}. Clearly, if $\{b_1=1_B, b_2,\cdots, b_n\}$ is a $k$-basis for the algebra $B$, then the projection $\bigoplus_{j=1}^nA\otimes_kb_j\ra A\otimes_k1_B$ shows that the injective map $A\ra A\otimes_kB$ is  left-split. Thus Corollary \ref{dmfortensor} follows from Theorem \ref{thmdm1}(2). $\square$

In general, the extension $A\ra A\otimes_kB$ does not have to be a Frobenius extension of $A$. For example, if we take $A :=k[x]/(x^2)$ and $B :=k(1\ra 2)$, the path algebra of $1\ra 2$ over $k$, then we can verify that $A\otimes_kB$ is not a Frobenius algebra, and therefore the extension is not a Frobenius extension since the algebra $A$ is Frobenius and since Frobenius extensions of algebras preserve Frobenius algebras (see \cite{Sugano, Nakayama2}).

\begin{Rem} \label{rmk}{\rm (1)
There is a more general result on dominant dimensions of tensor products of algebras by M\"uller in \cite[Lemma 6]{Muller}: For any finite-dimensional $k$-algebras $A$ and $B$ over a field $k$, $\dm(A\otimes_kB)$ = $min\{\dm(A),$ $\dm(B)\}$. Thus

$\quad$(i)  $\dm(A[x]/(x^n))=\dm(A)$ for all $n\ge 1$, where $A[x]$ is the polynomial algebra over $A$.

$\quad$(ii)  Let $M_n(A)$ be the $n\times n$ matrix algebra over $A$ and $T_n(A)$ be  the lower triangular $n\times n$ matrix subalgebra of $M_n(A)$. Then $\dm(T_n(A))=\min\{1, \dm(A)\}\le 1$ for $n>1$. Clearly, $M_n(A)$ is projective as a left or right $T_n(A)$-module. This shows that the condition ``left-split" in Lemma \ref{splitext} cannot be dropped.

(2) There are examples of Frobenius extensions $B\subseteq A$ such that $\dm(B)< \dm(A)$. For instance, we take Example 7.1 in \cite{Morita}. Let $A:=M_4(k)$, the $4\times 4$ matrix algebra over $k$, and $B$ be the subalgebra generated $k$-linearly by
$$ e_{11}+e_{44}, e_{22}+e_{33}, e_{21}, e_{31}, e_{41}, e_{42},e_{43}.$$
Then $B\subseteq A$ is a Frobenius extension. In this case, $\dm(B)=1 < \dm(A)=\infty$. Thus the  inequalities in Theorem \ref{thmdm1}(2) cannot be improved as equalities in general. Note that the Fronenius part of $B$ is the subalgebra $D$ generated by $e_{11}+e_{44}$ and $e_{41}$, and the Frobenius part of $A$ is $A$. Thus $D\subseteq A$ is a Frobenius extension by Theorem \ref{frobpart}.

(3) Suppose that $A, B, C$ and $D$ are $k$-algebras over a field $k$. If $_BM_A$ and $_DN_C$ are Frobenius bimodules, then $M\otimes_kN$ is a Frobenius $(B\otimes_kD)$-$(A\otimes_kC)$-bimodule. This can be deduced from \cite[Theorem 3.1, p. 209]{CE}, we leave the details to the reader.
So, if $B_i\subseteq A_i$ is a Frobenius extension of $k$-algebras for $i=1,2$, then $\dm(B_1\otimes_kB_2)\le \dm(A_1\otimes_kA_2)$ by Theorem \ref{thmdm1}(2).

(4) If $_BM_A$ and $_CL_B$ are Frobenius bimodules, then $_CL\otimes_BM_A$ is a Frobenius $C$-$A$-bimodule. This follows from the fact that $\Hom_A(_AX_C,\Hom_B(_DY_B,{}_AZ_B))\simeq \Hom_B(_DY_B,\Hom_A(_AX_C,{}_AZ_B))$ as $C$-$D$-bimodules for any $A$-$C$-bimodule $X$, $D$-$B$-bimodule $Y$ and $A$-$B$-bimodule $Z$. Thus the composition of Frobenius extensions is again a Frobenius extension.
}
\end{Rem}

Let $e=e^2$ be an idempotent element in a ring $R$. Suppose that $eR$ is a finitely generated projective $eRe$-module and the map: $Re\ra \Hom_{eRe}(eR,eRe), xe\mapsto (ey \mapsto eyxe)$ for $x,y\in R$, is bijective. Let $S:=eRe$. Then $eR$ is a Frobenius $S$-$R$-bimodule with $\add(eR)=\add(eRe)$. Thus, by Theorem \ref{thmdm1}, we have

\begin{Koro}
Let $A$ be a finite-dimensional algebra with an idempotent element $e\in A$. If $eA$ is a projective $eAe$-module and $\Hom_{eAe}(eA, eAe)\simeq Ae$ as $A$-modules, then $\dm(A)\le\dm(eAe)\le \dm(_AAe).$
\end{Koro}

For representation dimension, we conjectured in \cite[Section 7]{X-relative} that the representation dimension of $eAe$ is less than or equal to the one of $A$ for any idempotent element $e\in A$.

Recall that two algebras $A$ and $B$ are said to be \emph{stably equivalent of Morita type} if there are bimodules $_AU_B$ and $_BV_A$ such that they are projective as one-sided modules and have the properties: $_A U\otimes_BV_A\simeq {}_AA_A\oplus P$ and $_BV\otimes_AU_B\simeq {}_BB_B\oplus Q$ as bimodules, where $P$ and $Q$ are projective bimodules over $A$ and $B$, respectively.

In addition, if $A/\rad(A)$ and $B/\rad(B)$ are separable, then the bimodule $U$ and $V$ are Frobenius bimodules by a result in \cite{DM-V}. Since both $_AU$ and $_BV$ are generators, $\dm(A)=\dm(B)$. In fact, if $U$ and $V$ are indecomposable as bimodules, then $U$ and $V$ are still Frobenius bimodules even without the separability condition on quotients by
\cite[Lemma 4.2]{LX3} and \cite[Lamma 4.1]{CPX}. Thus we re-obtain the known fact that
the dominant dimensions of finite-dimensional algebras over a field are preserved by stable equivalences of Morita type. Note that stable equivalences of Morita type are not preserved by tensor products (see \cite{LZZ}). Nevertheless, we have

\begin{Koro} Suppose that $A$, $B$, $C$ and $D$ are finite-dimensional $k$-algebras over a field. If $A$ and $B$ as well as $C$ and $D$ are stably equivalent of Morita type, then $\dm(A\otimes_k C)=\dm(B\otimes_kD)$.
\label{moritatype}
\end{Koro}

{\it Proof.} Since $A$ and $B$ are stably equivalent of Morita type, $\dm(A)=\dm(B)$.
In general, we cannot get a stable equivalence of Morita type between $A\otimes_kC$ and $B\otimes_kC$ even for a symmetric $k$-algebra $C$ (see \cite{LZZ}). However, together with Remark \ref{rmk}(1), we have $$\dm(A\otimes_kC)=min\{\dm(A), \dm(C)\}=\min\{\dm(B),\dm(D)\} =\dm(B\otimes_kD). \square$$

\medskip
Next, we consider several cases of left-split Frobenius extensions.

Let $A$ be a $k$-algebra over a field $k$, and let $G$ be a finite group together with a group homomorphism
from $G$ to Aut$(A)$, the group of automorphisms of the $k$-algebra $A$. Let $kG$ be the
group algebra of $G$ over $k$. Then the skew group algebra $A\#_kG$ of $A$ by $G$
over $k$ has the underlying $k$-space $A\otimes_kkG$ with the multiplication given by
$$(a\otimes g)(b\otimes h) := a(b)g\otimes gh \quad \mbox{ for } \; a, b\in A; g, h\in G.$$
where $(b)g$ denotes the image of $b$ under $g$.

As a corollary of Theorem \ref{thmdm1}, we re-obtain \cite[Theorem 1.3(c)(ii)]{RR}.

\begin{Koro}(\cite{RR}) Let $A$ be a finite-dimensional $k$-algebra, and let $A\#_kG$ be the skew group algebra of $A$ by $G$
with $G$ a finite group. If $|G|$ is invertible in $A$, then $\dm(A\#_kG)$ = $\dm(A)$.
\end{Koro}

{\it Proof.} Note that under the assumption, the canonical inclusion $A\subseteq A\#_kG$ is a left-split Frobenius extension. $\square$

\medskip
Given an $A$-$A$-bimodule $M$, one defines the trivial extension $A\ltimes M$ of $R$ by $M$ as an algebra with the underlying $k$-space $A\oplus M$  and the multiplication: $(a,m)(b,n)=(ab, an+mb)$ for $a, b\in A; m,n\in M$. By Theorem  \ref{thmdm1}, we have

\medskip
\begin{Koro}If $M$ is an $A$-$A$-bimodule such that both $_AM$ and $M_A$ are projective, then $\dm(A\ltimes M)\le \dm(A)$. In particular, if $M\simeq Ae$ as bimodules for a central idempotent element $e\in A$, then $\dm(A\ltimes M)=\dm(A)$.\end{Koro}

{\it Proof.} The first statement follows immediately from Lemma \ref{splitext} since the canonical inclusion of $A$ into $A\ltimes M$ is left-split, while the second one is a consequence of Theorem \ref{thmdm1} because the canonical inclusion is a Frobenius extension by \cite[Theorem]{Kitamura} which says that $A\ltimes M$ is a Frobenius extension of $A$ if and only if there is a central idempotent element $e\in A$ such that $M\simeq Ae$ as $A$-$A$-bimodules. $\square$

\medskip
In connection with the Jones basic construction, there is a class of extensions, called Markov extensions.
The bulk of this material can be found in \cite[Chapter 3]{Kadison}. We just recall that an extension $B\subseteq A$ of $k$-algebras is called a \emph{Markov extension} if it is a Frobenius extension with a Frobenius system $(E,x_i,y_i)$, where $x_i, y_i\in A$ and $E:A\ra B$ is a $B$-$B$-bimodule homomorphism, and with a trace map $T: B\ra k$, such that $(1)E=1, \sum_i x_iy_i$ is nonzero element in $k$ and $E\circ T$ is a trace map on $A$. Markov extensions are left-split Frobenius extensions. Thus we have

\begin{Koro} If $B\subseteq A$ is a Markov extension of finite-dimensional algebras, then $\dm(A)=\dm(B)$.
\label{markovext}
\end{Koro}

\section{Flat-dominant dimensions of universal enveloping algebras of Lie algebras}

Following Morita \cite{moritaqf3}, a ring $R$ is called a \emph{left
(resp. right) $QF$-$3$} ring if $\fdd(_RR)\ge 1$ (resp. $\fdd(R_R)$ $\ge 1$). By a $QF$-$3$ ring we mean a left and right $QF$-$3$ ring. Note that $QF$-$3$ rings possess remarkable properties, for instance, for a noetherian ring $R$, it is left $QF$-$3$
if and only if the functor $\Hom_R(-,\Hom_R(-,R )): R\Modcat\ra R\Modcat$ preserves monomorphisms (see \cite{Hoshino, T}).

In this section, we show that the quantised enveloping algebras of semisimple Lie algebras are $QF$-$3$ rings.

Let $n$ be a natural number. Recall that a commutative noetherian ring $R$ is said to be \emph{$n$-Gorenstein} if the injective dimensions of $_RR$ is equal to $n$.

To apply our result to various algebras in other fields, we need a couple of lemmas.

\begin{Lem} Suppose that $S$ is a commutative noetherian ring and $\lambda: S\ra R$ is a localization of $S$ (with respect to a multiplicative subset of $S$). Then

$(1)$ $\fdd(S)\le \fdd(R)$.

$(2)$ Let $S$ be an $n$-Gorenstein ring with a minimal injective resolution of $_SS: 0\ra S\ra I_0\ra \cdots\ra I_{n-1}\ra I_n\ra 0$. Then the injective dimension of $_RR$ is at most $n$, and it is $n$ if and only if $R\otimes_SI_j\ne 0$ for all $I_j$ with $0\le j\le n-1$.
\label{local}
\end{Lem}

{\it Proof.} (1) Since the ring $S$ is noetherian, the localization $R$ of $S$ is again noetherian. For noetherian rings, localizations preserve injective modules. Thus $\fdd(S)\le \fdd(R)$.

(2) follows from the result (see \cite[Theorem 1]{Iv}):
Let $A$ be a commutative noetherian ring, and let $F$ be a flat $A$-module. If $f: M\ra N$ is an essential
monomorphism, then $1_F\otimes f: F\otimes_A M\ra F\otimes_AN$ is an essential monomorphism.
$\square$

\begin{Lem} \label{gorensteinfdd} If $R$ is a commutative $n$-Gorenstein ring with $n\ge 1$, then
$\fdd(R)=1$. In particular, $\fdd(R[X_1, \cdots, X_m])=1$ and $\fdd(R[X_1^{\pm}, \cdots, X_m^{\pm}])=1$.
\end{Lem}

{\it Proof.} We have a minimal injective resolution of the $n$-Gorenstein ring $R$:
$$0\lra R\lra I_0\lra \cdots\lra I_n\lra 0, \; \; I_i =\bigoplus_{p,\; ht(p)=i}E(R/p) \; \mbox{ for } \; 0\le i\le n,$$
where $p$ runs over prime ideals in $R$, $ht(p)$ stands for the height of $p$, and $E(R/p)$ is the injective hull of the $R$-module $R/p$. By \cite{Xu}, $E(R/p)$ is flat if and only if the height of $p$ is $0$. Thus $\fdd(R)=1$. Clearly, $R[X_1, \cdots, X_m]$ is a $(n+m)$-Gorenstein ring, and therefore $\fdd(R[X_1, \cdots, X_m])=1$. Since $R[X_1^{\pm}, \cdots, X_m^{\pm}]$ is a localization of
$R[X_1, \cdots, X_m]$, it is  an $(n+m)$-Gorenstein ring and has $\fdd(R[X_1^{\pm}, \cdots, X_m^{\pm}])=1$ by Lemma \ref{local}(2). $\square$

\begin{Koro}
Let $\mathfrak{g}$ be a finite-dimensional semisimple Lie algebra of rank $r$ over a field $k$ and $U(\mathfrak{g})$ the universal enveloping algebra of $\mathfrak{g}$.
If $r\ge 1$, then $\fdd(U(\mathfrak{g}))=1$. Thus $U(\mathfrak{g})$ is a $QF$-$3$ ring
\end{Koro}

{\it Proof.} By \cite{BrownGordonStroppel}, $U(\mathfrak{g})$ is a free Freobenius extension of its center $\mathcal{Z}$, while the center $\mathcal{Z}$ of $U(\mathfrak{g})$ is isomorphic to the polynomial algebra $k[X_1,\cdots, X_r]$ in $r$ variables (see, for example,  \cite[Theorem 8]{XavierBekaert}). Thus, by
Theorem \ref{thmdm1}(1), $\fdd(U(\mathfrak{g}))=\fdd(k[X_1,\cdots, X_r])=1.$ Since $\mathcal{Z}$ is a commutative ring, the right flat-dominant dimension of $U(\mathfrak{g})$ can be calculated similarly. $\square$

\medskip
The above result can be extended to quantised enveloping algebras of Lie algebras.
We recall briefly a few facts on the quantised enveloping algebra of a semisimple Lie algebra from \cite[III 6.1-6.2]{browngoodearl}. Assume that $\mathfrak{g}$ is a finite-dimensional semisimple
complex Lie algebra and that $\epsilon$ is a primitive $\ell$-th root of unity in the field $k$ of
characteristic $0$, where
$$\ell>3 \mbox{ is odd, and prime to } 3 \mbox{ if } \mathfrak{g} \mbox{ contains a factor of type } G_2.$$
Let $\check{U}_{\epsilon}(\mathfrak{g})$ be the quantised enveloping algebra of $\mathfrak{g}$ over $k$, and $Z_0$ be the sub-Hopf centre of $\check{U}_{\epsilon}(\mathfrak{g})$. Then $\check{U}_{\epsilon}(\mathfrak{g})$
is a free Frobenius extension of $Z_0$ of rank $\ell^{\mbox{dim}(\mathfrak{g})}$ (see \cite{BrownGordonStroppel}) and $Z_0$ is isomorphic to a polynomial algebra over $k$ in $2N+n$ variables with $n$ variable inverted, where $N$ is the number of positive roots of $\mathfrak{g}$ (see \cite[III, 6.2]{browngoodearl}).
Thus $\fdd(Z_0)=1$ by Lemmas \ref{local} and \ref{gorensteinfdd}. Note that $\check{U}_{\epsilon}(\mathfrak{g})$ is a noetherian ring. Hence the following corollary holds by Theorem \ref{thmdm1}(1).

\begin{Koro} $\check{U}_{\epsilon}(\mathfrak{g})$ is a $QF$-$3$ ring. In fact, $\fdd(\check{U}_{\epsilon}(\mathfrak{g}))=1$.
\label{quantised}
\end{Koro}

\medskip
{\bf Acknowledgements.} The research work of the author was partially supported by the Beijing Natural Science Foundation (1192004) and the National Natural Science Foundation of China. The author is grateful to the Alexander von Humboldt Foundation for supporting an academic visit to the University of Stuttgart in the summer of 2018, where part of the work was revised.

\medskip
{\footnotesize
Changchang Xi, School of Mathematical Sciences, Capital Normal University, Beijing 100048, China; and School of Mathematics and Information Science,
Henan Normal University, Xinxiang 453007, Henan, China

{\tt Email: xicc@cnu.edu.cn}
\end{document}